\documentclass[a4]{article}
\usepackage{graphicx}

\def\be{\begin{equation}}
\def\ee{\end{equation}}
\def\a{\alpha}
\def\b{\beta}

\def\x{\mathbf{x}}

\begin{document}

\title{Firefly Algorithm, Stochastic Test Functions and Design Optimisation}

\author{Xin-She Yang \\
Department of Engineering, University of Cambridge, \\
Trumpington Street, Cambridge CB2 1PZ, UK \\
Email: xy227@cam.ac.uk
}

\maketitle


\begin{abstract}
Modern optimisation algorithms are often metaheuristic, and they are
very promising in solving NP-hard optimization problems. In this paper, we show how
to use the recently developed Firefly Algorithm to solve nonlinear design problems.
For the standard pressure vessel design optimisation, the optimal solution
found by FA is far better than the best solution obtained
previously in literature. In addition, we also propose
a few new test functions with either singularity or stochastic components but with
known global optimality, and thus they can be used to validate new optimisation algorithms.
Possible topics for further research are also discussed.
\end{abstract}

{\it To cite this paper as follows: Yang, X. S., (2010)
`Firefly Algorithm, Stochastic Test Functions and Design Optimisation',
{\it Int. J. Bio-Inspired Computation}, Vol.~2, No.~2, pp.78--84.}

\section{Introduction}

Most optimization problems in engineering are nonlinear with many constraints.
Consequently, to find optimal solutions to such nonlinear problems requires
efficient optimisation algorithms (Deb 1995, Baeck et al 1997, Yang 2005).
In general, optimisation algorithms can be
classified into two main categories: deterministic and stochastic.
Deterministic algorithms such as hill-climbing will produce the same set of
solutions if the iterations start with the same initial guess. On the other hand,
stochastic algorithms often produce different solutions even with
the same initial starting point. However, the final results, though slightly different,
will usually converge to the same optimal solutions within a given accuracy.

Deterministic algorithms are almost all local search algorithms, and they are quite
efficient in finding local optima. However, there is a risk
for the algorithms to be trapped at local optima, while the global optima are
out of reach. A common practice is to introduce some stochastic component to
an algorithm so that it becomes possible to jump out of such
locality. In this case, algorithms become stochastic.

Stochastic algorithms often have a deterministic component and
a random component. The stochastic component can take many forms such as
simple randomization by randomly sampling the search space or by random walks.
Most stochastic algorithms can be considered as metaheuristic, and good examples
are genetic algorithms (GA) (Holland 1976, Goldberg 1989) and particle swarm optimisation (PSO)
(Kennedy and Eberhart 1995, Kennedy et al 2001).
Many modern metaheuristic algorithms
were developed based on the swarm intelligence in nature (Kennedy and Eberhart 1995,
Dorigo and St\"utzle 2004).
New modern metaheuristic algorithms are being developed and begin
to show their power and efficiency. For example, the Firefly Algorithm
developed by the author shows its superiority over some  traditional algorithms
(Yang 2009, Yang 2009, Lukasik and Zak 2009).

The paper is organized as follows: we will first briefly outline the
main idea of the Firefly Algorithm in Section 2, and  we then describe
a few new test functions with singularity and/or randomness in Section 3.
In Section 4, we will use FA to find the optimal solution of a pressure vessel design problem.
Finally, we will discuss the topics for further studies.

\section{Firefly Algorithm and its Implementation}

\subsection{Firefly Algorithm}

The Firefly Algorithm was developed by the author (Yang 2008, Yang 2009),
and it was based on the idealized behaviour of the
flashing characteristics of fireflies. For simplicity, we
can idealize these flashing characteristics as the following
three rules
\begin{itemize}

\item all fireflies are unisex so that one firefly is attracted to other fireflies
regardless of their sex;

\item Attractiveness is proportional to their brightness,
thus for any two flashing fireflies, the less brighter one will move towards the brighter
one. The attractiveness is proportional to the brightness and they both
decrease as their distance increases. If
no one is brighter than a particular firefly, it moves randomly;

\item The brightness or light intensity of a firefly is affected or  determined by the landscape of the
objective function to be optimised.

\end{itemize}

For a  maximization problem, the brightness can simply be proportional
to the objective function. Other forms of brightness can be defined in a similar
way to the fitness function in genetic algorithms or the bacterial
foraging algorithm (BFA) (Gazi and Passino 2004).

In the FA, there are two important issues: the variation of
light intensity and formulation of the attractiveness.
For simplicity, we can always assume that the attractiveness
of a firefly is determined by its brightness or light intensity
which in turn is associated with the encoded objective function.
In the simplest case for maximum optimization problems,
the brightness $I$ of a firefly at a particular location $\x$ can be chosen
as $I(\x) \propto f(\x)$. However, the attractiveness $\b$ is relative, it should be
seen in the eyes of the beholder or judged by the other fireflies. Thus, it should
vary with the distance $r_{ij}$ between firefly $i$ and firefly $j$.
As light intensity decreases with the distance from its source, and light is also
absorbed in the media,
so we should allow the attractiveness to vary with the degree of absorption.

In the simplest form, the light intensity $I(r)$ varies with the distance $r$
monotonically and exponentially. That is
\be I=I_0 e^{-\gamma r}, \ee
where $I_0$ is the original light intensity
and $\gamma$ is the light absorption coefficient.
As a firefly's attractiveness is proportional to
the light intensity seen by adjacent fireflies, we can now define
the attractiveness $\b$ of a firefly by
\be \b = \b_0 e^{-\gamma r^2}, \label{att-equ-100} \ee
where $\b_0$ is the attractiveness at $r=0$. It is worth pointing out
that the exponent $\gamma r$ can be replaced by other functions
such as $\gamma r^m$ when $m>0$. Schematically,
the Firefly Algorithm (FA) can be summarised as the pseudo code.

\begin{center}
\begin{minipage}{0.77\textwidth}
{\sf
Firefly Algorithm
\hrule
\vspace{5pt}
\indent Objective function $f(\x), \quad \x=(x_1, ..., x_d)^T$ \\
\indent Initialize a population of fireflies $\x_i \; (i=1,2,...,n)$ \\
\indent Define light absorption coefficient $\gamma$ \\
\indent {\bf while} ($t<$MaxGeneration) \\
\indent {\bf for} $i=1:n$ all $n$ fireflies \\
\indent \quad {\bf for} $j=1:i$ all $n$ fireflies \\
\indent \qquad Light intensity $I_i$ at $\x_i$ is determined by $f(\x_i)$  \\
\indent \qquad  {\bf if} ($I_j>I_i$) \\
\indent \qquad Move firefly $i$ towards $j$ in all $d$ dimensions \\
\indent \qquad {\bf end if} \\
\indent \qquad Attractiveness varies with distance $r$ via $\exp[-\gamma r]$ \\
\indent \qquad Evaluate new solutions and update light intensity \\
\indent \quad {\bf end for }$j$ \\
\indent {\bf end for }$i$ \\
\indent Rank the fireflies and find the current best  \\
\indent {\bf end while} \\
\indent Postprocess results and visualization }
\hrule
\end{minipage}
\end{center}

The distance between any two fireflies $i$ and $j$ at $\x_i$ and $\x_j$
can be the Cartesian distance  $r_{ij}=||\x_i-\x_j||_2$ or the $\ell_2$-norm.
For other applications such as scheduling, the distance can be time delay or
any suitable forms, not necessarily the Cartesian distance.

The movement of a firefly $i$ is attracted to another more attractive (brighter)
firefly $j$ is determined by
\be \x_i =\x_i + \b_0 e^{-\gamma r^2_{ij}} (\x_j-\x_i) + \a \mbox{\boldmath $\epsilon$}_i, \ee
where the second term is due to the attraction, while the third term
is randomization with the vector of random variables $\mbox{\boldmath $\epsilon$}_i$ being  drawn from a Gaussian distribution.

For most cases in our implementation,
we can take $\b_0=1$, $\a \in [0,1]$, and $\gamma=1$.
In addition, if the scales vary significantly in different dimensions such as $-10^5$ to $10^5$ in
one dimension while, say, $-10^{-3}$ to $10^3$ along others, it is a good idea to
replace $\a$ by $\a S_k$ where the scaling parameters $S_k (k=1,...,d)$ in
the $d$ dimensions should be determined by the actual scales of the problem of interest.

In essence, the parameter $\gamma$ characterizes the variation of the attractiveness,
and its value is crucially important in determining the speed of the convergence
and how the FA algorithm behaves. In theory, $\gamma \in [0,\infty)$, but in
practice, $\gamma=O(1)$ is determined by the characteristic/mean length $S_k$ of the
system to be optimized. In one extreme when $\gamma \rightarrow 0$,
the attractiveness is constant $\b=\b_0$. This is equivalent to saying
that the light intensity does not decrease in an idealized sky.
Thus, a flashing firefly can be seen anywhere in the domain. Thus, a single (usually global)
optimum can easily be reached. This corresponds to a special case of particle
swarm optimization (PSO). In fact, if the inner loop
for $j$ is removed and $I_j$ is replaced by the current global best $\mathbf{g}_*$,
FA essentially becomes the standard PSO, and, subsequently,
the efficiency of this special case is the same as that of PSO.
On the other hand, if $\gamma \rightarrow \infty$, we have
$\b(r) \rightarrow \delta(r)$, which is a Dirac $\delta$-function.
This means that  the attractiveness is almost zero in the sight of
other fireflies, or the fireflies are short-sighted. This is equivalent
to the case where the fireflies fly in a very foggy region randomly. No other fireflies can be
seen, and each firefly roams in a completely random way. Therefore, this corresponds
to the completely random search method. So $\gamma$ partly controls how
the algorithm behaves. It is also possible to adjust $\gamma$ so that
multiple optima can be found at the same during iterations.

\subsection{Numerical Examples}

From the pseudo code, it is relatively straightforward to implement the Firefly Algorithm
using a popular programming language such as Matlab. We have tested it against more than
a dozen test functions such as the Ackley function
\[ f(\x)=-20 \exp\Big[-\frac{1}{5} \sqrt{\frac{1}{d} \sum_{i=1}^d x_i^2}\Big] \]
\be - \exp[\frac{1}{d} \sum_{i=1}^d \cos (2 \pi x_i)] +20 +e, \ee
which has a unique global minimum $f_*=0$ at $(0,0,...,0)$.
From a simple parameter studies, we concluded that,
in our simulations, we can use the following values of
parameters $\a=0.2$, $\gamma=1$, and $\b_0=1$.
\begin{figure}
 \centerline{\includegraphics[width=4in,height=3in]{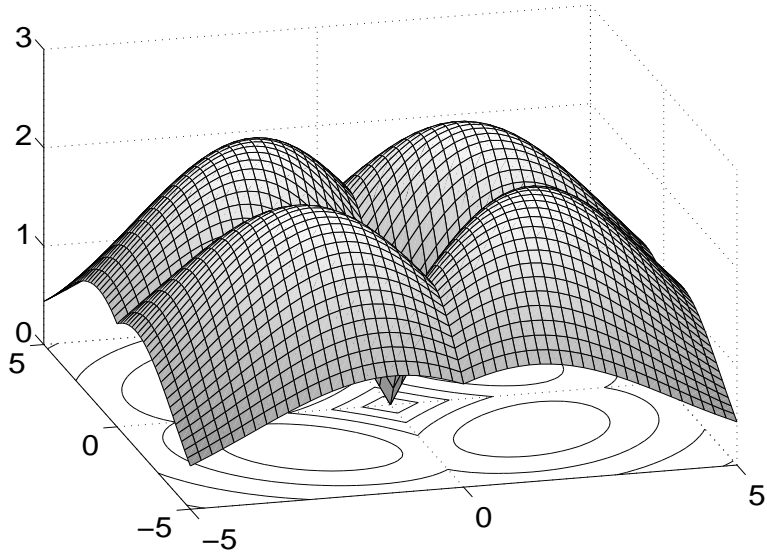} }
\caption{Four global maxima at $(\pm 1/2,\pm 1/2)$. \label{Yang-fig1} }
\end{figure}
As an example, we now use the FA to find the
global maxima of the following function
\be f(\x)=\Big( \sum_{i=1}^d |x_i| \Big) \exp \Big( -\sum_{i=1}^d x_i^2 \Big),
\label{fun-equ-50} \ee
with the domain $-10 \le x_i \le 10$ for all $(i=1,2,...,d)$
where $d$ is the number of dimensions. This function has
multiple global optima.  In the case of $d=2$, we have $4$ equal maxima
$f_* =1/\sqrt{e} \approx 0.6065$ at $(1/2,1/2)$, $(1/2,-1/2)$,
$(-1/2,1/2)$ and $(-1/2,-1/2)$ and a unique global minimum at $(0,0)$.

\begin{figure}
 \centerline{\includegraphics[height=3in,width=4in]{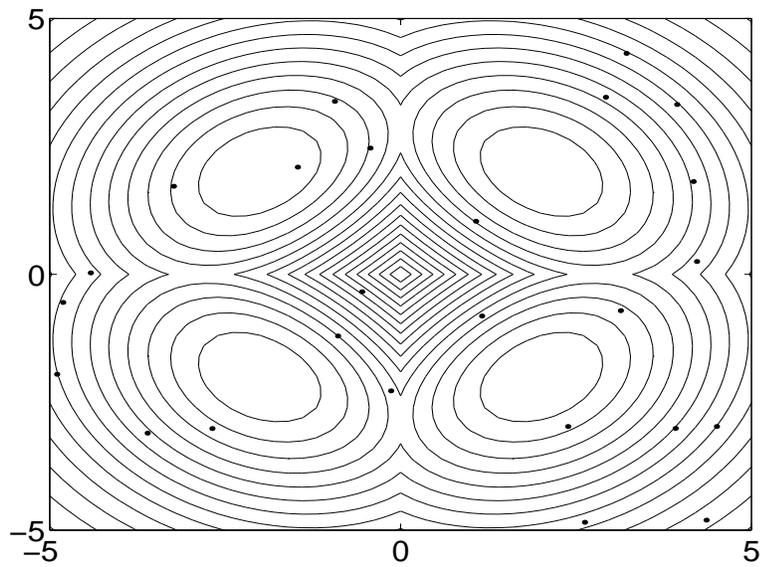}}
\vspace{-5pt}
\caption{Initial locations of 25 fireflies. \label{Yang-fig2a} }
\end{figure}

\begin{figure}
 \centerline{\includegraphics[height=3in,width=4in]{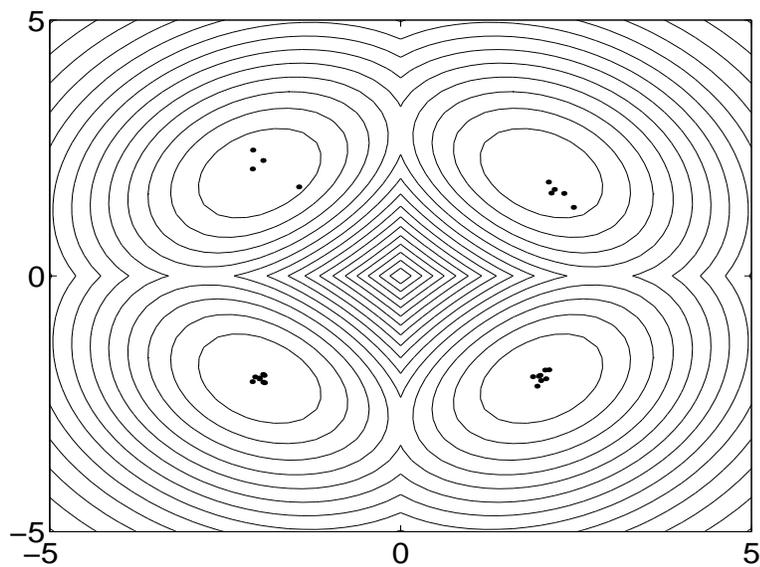}}
\vspace{-5pt}
\caption{Final locations after $20$ iterations. \label{Yang-fig2b} }
\end{figure}

The four peaks are shown in Fig. \ref{Yang-fig1}, and these global maxima
can be found using the implemented Firefly Algorithms after about 500 function
evaluations. This corresponds to $25$ fireflies evolving for $20$ generations or
iterations. The initial locations of 25 fireflies are shown Fig. \ref{Yang-fig2a}
and their final locations after 20 iterations are shown in Fig. \ref{Yang-fig2b}.
We can see that the Firefly Algorithm is very efficient. Recently studies
also confirmed its promising power in solving nonlinear
constrained optimization tasks (Yang 2009, Lukasik and Zak 2009).

\section{New Test Functions}

The literature about test functions is vast, often with different collections of
test functions for validating new optimisation algorithms. Test
functions such as Rosenbrock's banana function and Ackley's
function mentioned earlier are well-known in the optimisation literature.
Almost all these test functions are deterministic and smooth.
In the rest of
this paper, we first propose a few new test functions which have
some singularity and/or stochastic components. Some of the formulated
functions have stochastic components but their global optima are deterministic.
Then, we will use the Firefly Algorithm to find the optimal solutions of some
of these new functions.

\begin{figure}
 \centerline{\includegraphics[width=4in,height=3in]{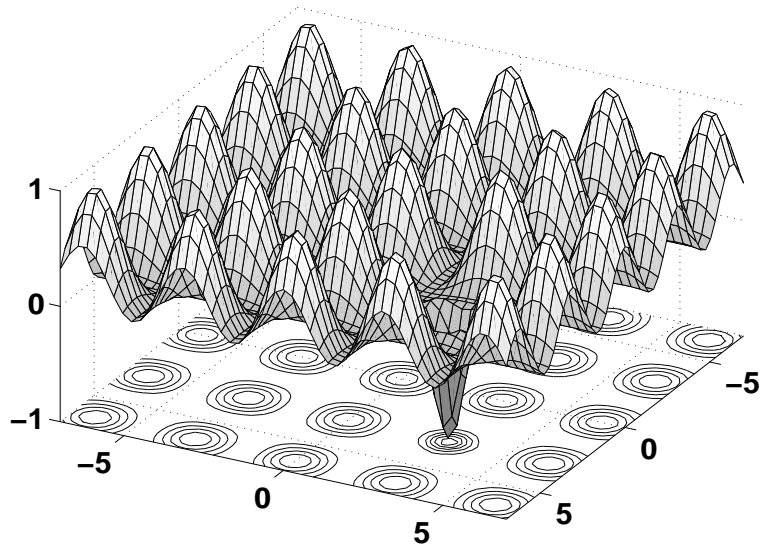}}
\caption{The standing wave function for two independent variables
with the global minimum $f_*=-1$ at $(\pi,\pi)$. \label{Yang-fig3} }
\end{figure}

\begin{figure}
 \centerline{\includegraphics[width=4in,height=3in]{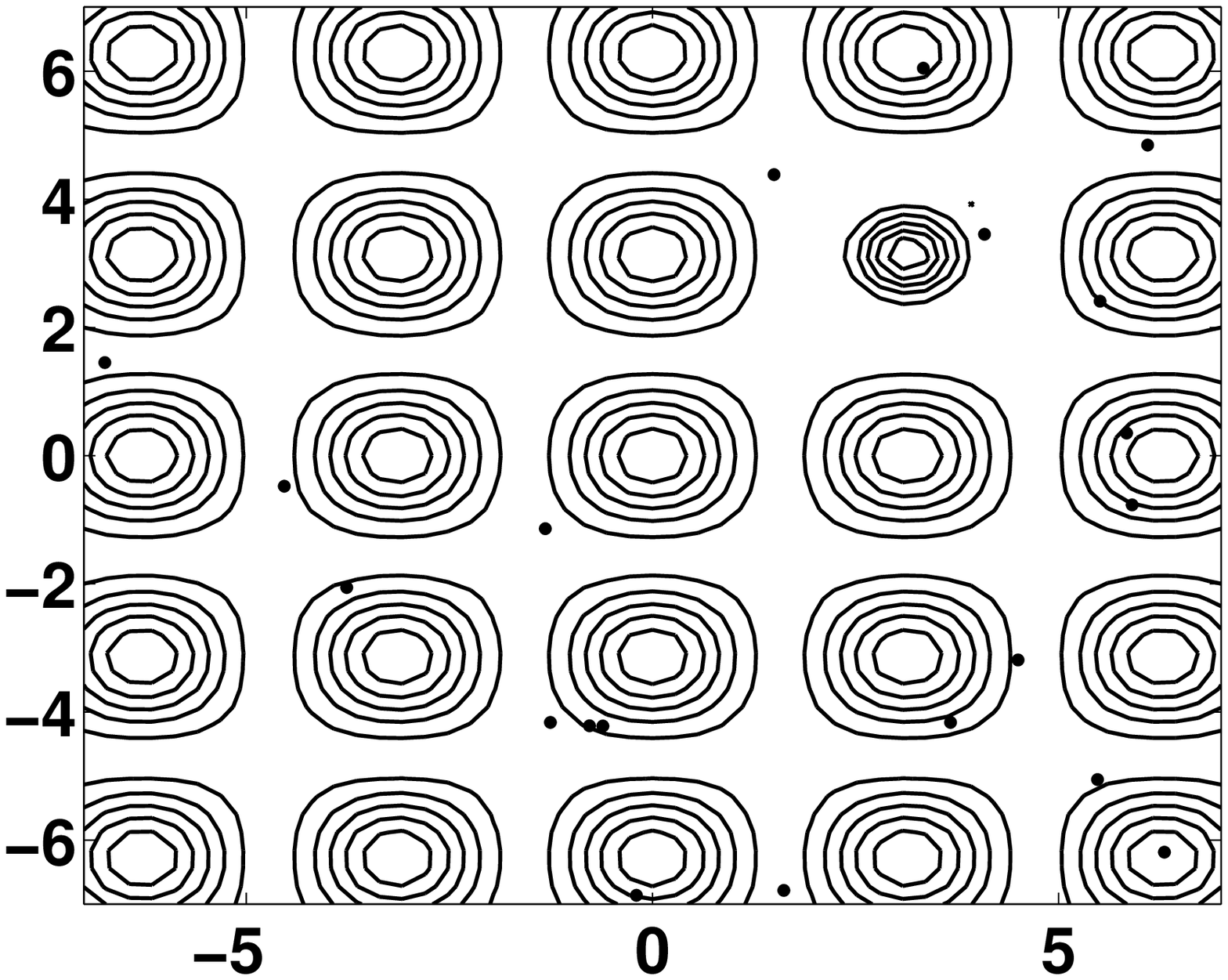} }

\caption{The initial locations of $20$ fireflies.  \label{Yang-fig4a} }
\end{figure}

\begin{figure}
 \centerline{\includegraphics[width=4in,height=3in]{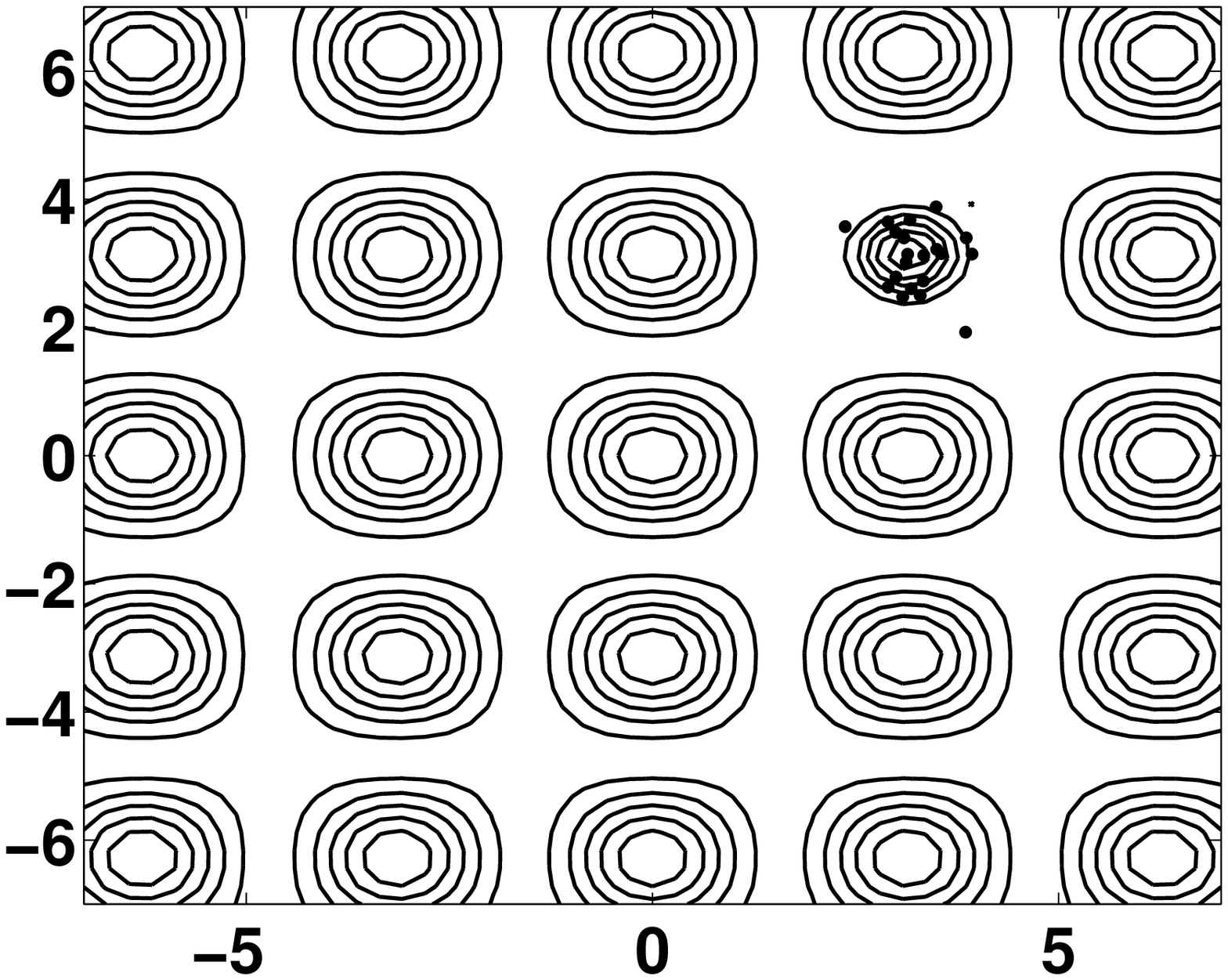} }

\caption{Final locations after $15$ iterations.  \label{Yang-fig4b} }
\end{figure}

The first test function we have designed is a multimodal nonlinear function
\[ f(\x)=\Big[ e^{-\sum_{i=1}^d (x_i/\b)^{2m}} - 2 e^{-\sum_{i=1}^d (x_i-\pi)^2} \Big] \]
\be \cdot \prod_{i=1}^d \cos^2 x_i, \quad m=5, \ee
which looks like a standing-wave function with a defect (see Fig. \ref{Yang-fig3}).
It has many local minima and the unique global minimum $f_*=-1$ at $\x_*=(\pi,\pi,...,\pi)$
for $\b=15$ within the domain $-20 \le x_i \le 20$ for $i=1,2,...,d$.
By using the Firefly Algorithm with $20$ fireflies, it is easy to
find the global minimum in just about $15$ iterations. The results are
shown in Fig. \ref{Yang-fig4a} and Fig. \ref{Yang-fig4b}.

\begin{figure}
 \centerline{\includegraphics[width=4in,height=3in]{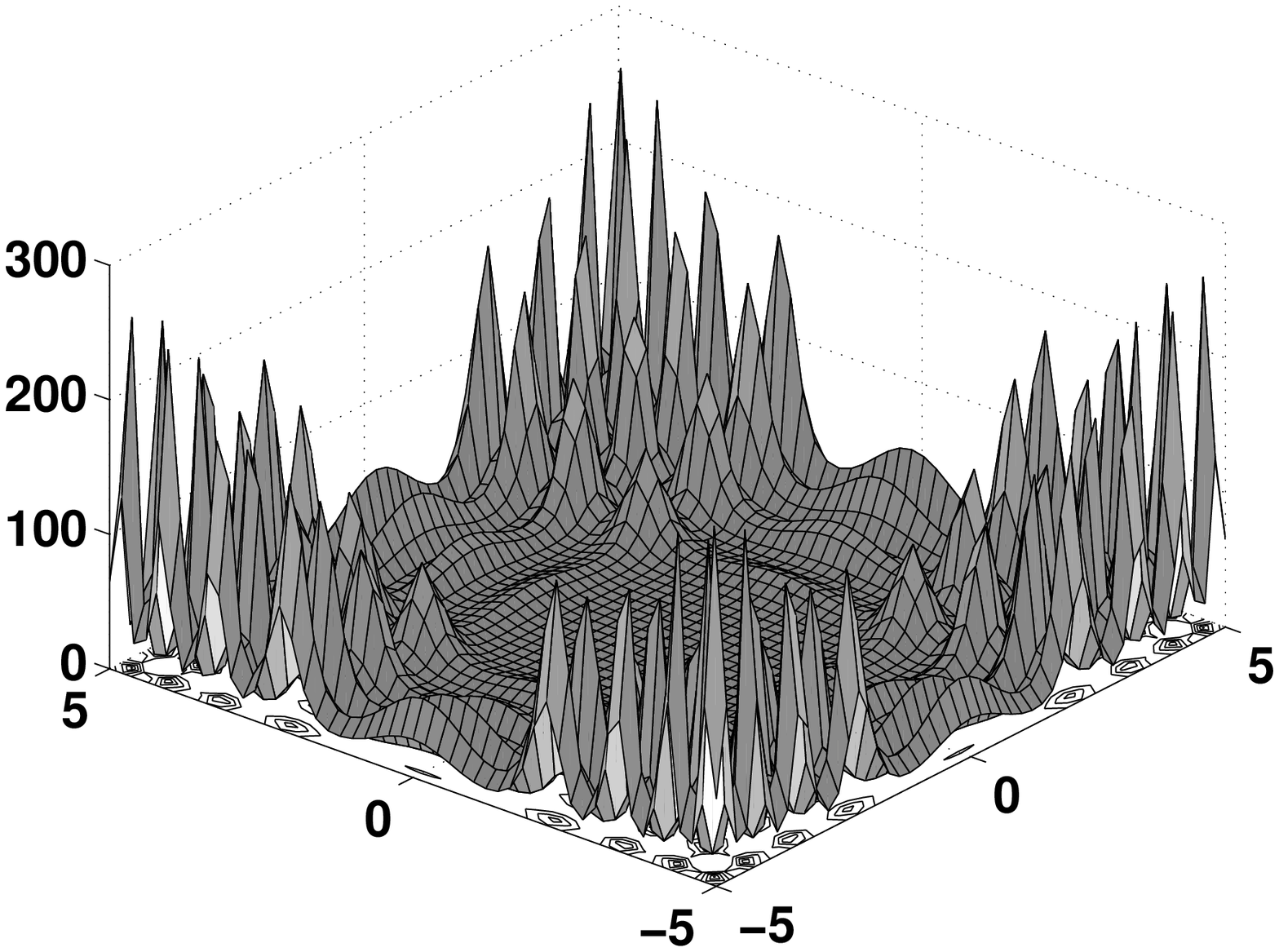} }
\caption{The landscape of function (\ref{fun-equ-100}). \label{Yang-fig5a} }
\end{figure}

\begin{figure}
 \centerline{\includegraphics[width=4in,height=3in]{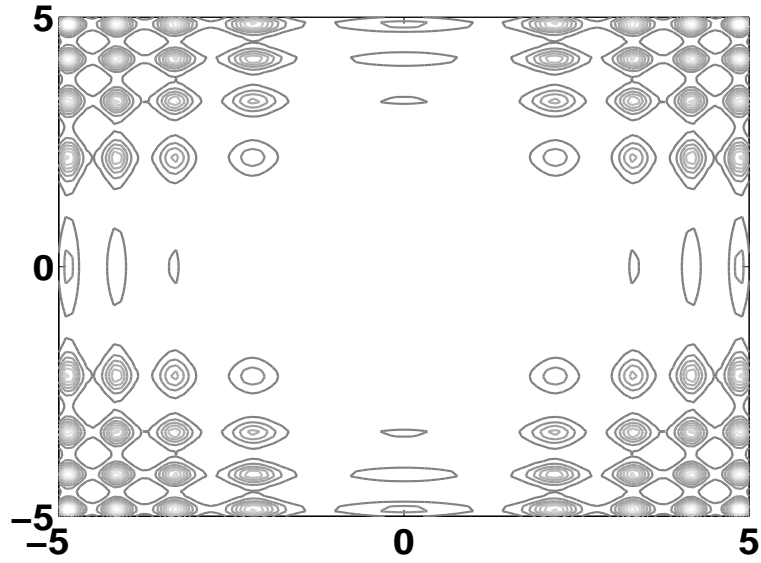} }
\caption{Contour of function (\ref{fun-equ-100}). \label{Yang-fig5b} }
\end{figure}

As most test functions are smooth, the next function we have formulated
is also multimodal but it has a singularity
\be f(\x) =\Big( \sum_{i=1}^d |x_i| \Big)  \cdot \exp \Big[ -\sum_{i=1}^d \sin (x_i^2) \Big],
\label{fun-equ-100} \ee
which has a unique global minimum $f_*=0$ at $\x_*=(0,0,...,0)$
in the domain $-2 \pi \le x_i \le 2 \pi$ where $i=1,2,...,d$.
At a first look,  this function has some similarity with function (\ref{fun-equ-50}) discussed earlier.
However, this function is not smooth,  and its derivatives are not well defined at the optimum $(0,0,...,0)$.
The landscape of this forest-like function is shown in Fig. \ref{Yang-fig5a}
and its 2D contour is displayed in Fig. \ref{Yang-fig5b}.

Almost all existing test functions are deterministic. Now let us
design a test function with stochastic components
\[ f(x,y)=-5 e^{-\b  [(x-\pi)^2+(y-\pi)^2]} \]
\be -\sum_{j=1}^K \sum_{i=1}^K \epsilon_{ij} e^{-\a  [(x-i)^2 + (y-j)^2]}, \ee
where $\a,\b>0$ are scaling parameters, which can often be taken as $\a=\b=1$.
Here the random variables $\epsilon_{ij} \;(i,j=1,...,K)$ obey
a uniform distribution $\epsilon_{ij} \sim $ Unif[0,1].
The domain is $0 \le x,y \le K$ and $K=10$. This function has
$K^2$ local valleys at grid locations and the  fixed global
minimum at $\x_*=(\pi,\pi)$. It is worth pointing that
the minimum $f_{\min}$ is random, rather than a fixed value;
it may vary from $-(K^2+5)$ to $-5$, depending $\a$ and $\b$ as well as the
random numbers drawn.

\begin{figure}
 \centerline{\includegraphics[width=4in,height=3in]{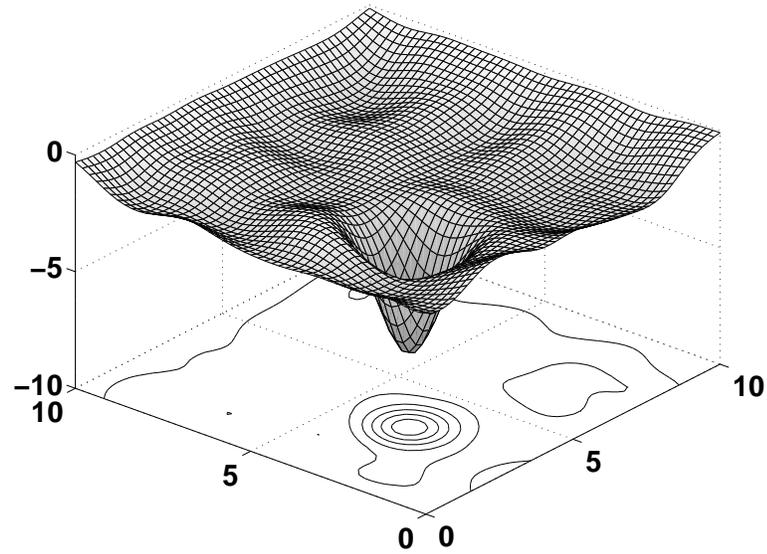}}
\caption{The 2D Stochastic function for $K=10$
with a unique global minimum at $(\pi,\pi)$, though
the value of this global minimum is somewhat stochastic. \label{Yang-fig6} }
\end{figure}

For stochastic test functions, most deterministic algorithms such as hill-climbing
would simply fail due to the fact that the landscape is constantly changing.
However, metaheuristic algorithms could still be robust in dealing with
such functions. The landscape of a realization of this stochastic function is shown in
Fig. \ref{Yang-fig6}.

\begin{figure}
 \centerline{\includegraphics[width=4in,height=3in]{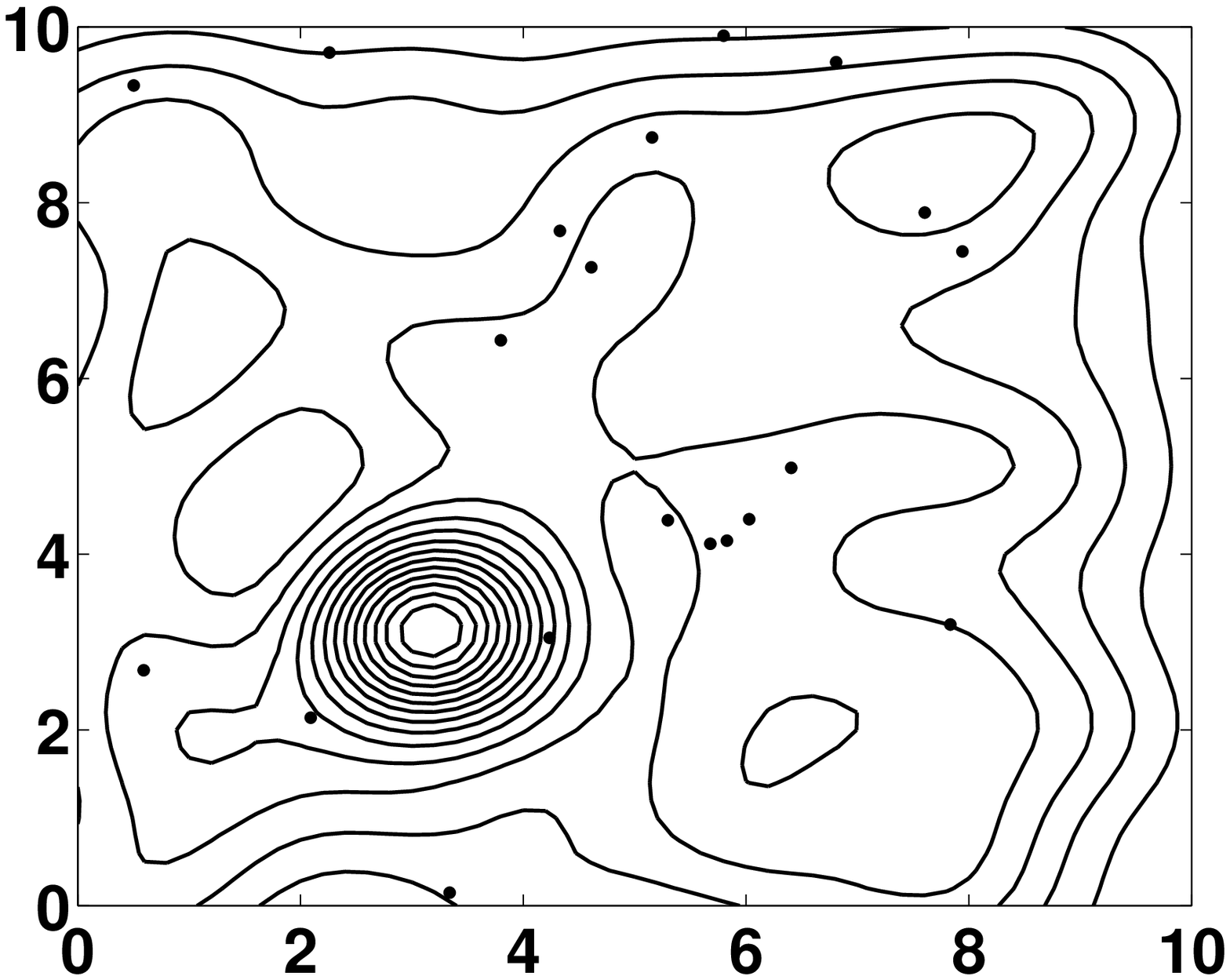} }
\caption{The initial locations of 20 fireflies. \label{Yang-fig7a} }
\end{figure}

\begin{figure}
 \centerline{ \includegraphics[width=4in,height=3in]{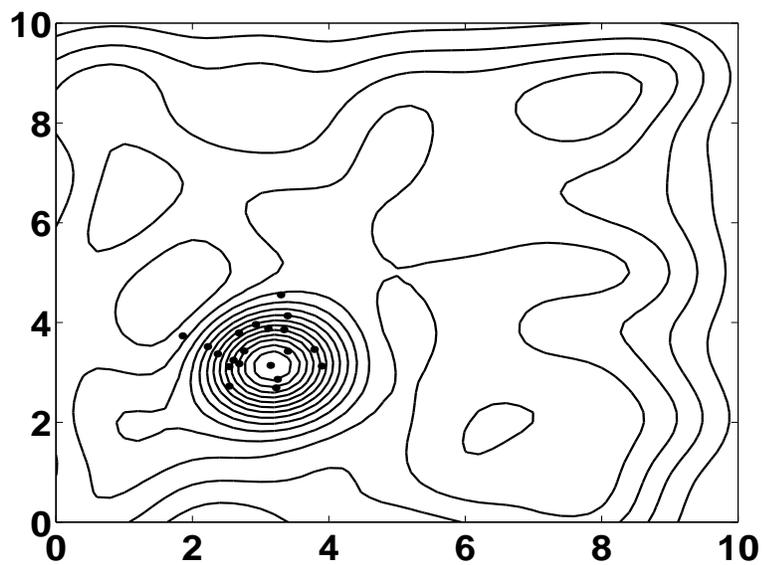} }
\caption{The final locations of 20 fireflies after $15$ iterations, converging
into $(\pi,\pi)$. \label{Yang-fig7b} }
\end{figure}

Using the Firefly Algorithm, we can find the
global minimum in about 15 iterations for $n=20$ fireflies. That is,
the total number of function evaluations is just 300. This is indeed
very efficient and robust. The initial locations of the fireflies are
shown in Fig. \ref{Yang-fig7a} and the final results are shown Fig.\ref{Yang-fig7b}.

Furthermore, we can also design a relative generic stochastic function
which is both stochastic and non-smooth
\be f(\x) = \sum_{i=1}^d \epsilon_i \; \Big|x_i\Big|^i, \quad -5 \le x_i \le 5, \ee
where $\epsilon_i \;(i=1,2,...,d)$ are random variables which are uniformly distributed
in $[0,1]$. That is, $\epsilon_i \sim $Unif$[0,1]$.
This function has the unique minimum $f_*=0$ at $\x_*=(0,0,...,0)$
which is also singular.

 We found that for most problems $n=15$ to $50$ would be sufficient.
For tougher problems, larger $n$ can be used, though excessively
large $n$ should not be used unless there is no better alternative,
as it is more computationally extensive for large $n$.

\section{Engineering Applications}

Now we can apply the Firefly Algorithm to carry out various design optimisation
tasks. In principle, any optimization problems that can be solved by genetic
algorithms and particle swarm optimisation can also be solved by the Firefly Algorithm.
For simplicity to demonstrate its effectiveness in real-world optimisation, we use the FA
to find the optimal solution of the standard but quite difficult pressure design
optimisation problem.

Pressure vessels are literally everywhere such as champagne bottles and gas tanks.
For a given volume and working pressure, the basic aim of designing
a cylindrical vessel is to minimize the total cost. Typically, the design variables are
the thickness $d_1$ of the head, the thickness $d_2$ of the body, the inner radius
$r$, and the length $L$ of the cylindrical section (Coello 2000, Cagnina et al 2008).
This is a well-known test problem for optimization  and it can be written as
\[ \textrm{minimize } f(\x) = 0.6224 d_1 r L + 1.7781 d_2 r^2 \]
\be + 3.1661 d_1^2 L + 19.84 d_1^2 r, \ee
subject to the following constraints
\be
\begin{array}{lll}
 g_1(\x) = -d_1 + 0.0193 r \le 0 \\
 g_2(\x) = -d_2 + 0.00954 r \le 0 \\
 g_3(\x) = - \pi r^2 L -\frac{4 \pi}{3} r^3 + 1296000 \le 0 \\
 g_4(\x) =L -240 \le 0.
\end{array}
\ee

The simple bounds are  \be 0.0625 \le d_1, d_2 \le 99 \times 0.0625, \ee
and \be 10.0 \le r, \quad L \le 200.0. \ee

Recently, Cagnina et al (2008) used an efficient particle swarm optimiser
to solve this problem and they found the best solution
\be f_* \approx 6059.714, \ee
at \be \x_* \approx (0.8125, \; 0.4375, \; 42.0984, \; 176.6366). \ee
This means the lowest price is about $\$6059.71$.

Using the Firefly Algorithm, we have found an even better solution
with $40$ fireflies after $20$ iterations, and we have obtained
\be \x_* \approx (0.7782, 0.3846, 40.3196, 200.0000)^T, \ee
with \be f_{\min} \approx 5885.33, \ee
which is significantly lower or cheaper than the solution
$f_* \approx 6059.714$ obtained by Cagnina et al (2008).

This clearly shows how efficient and effective the Firefly Algorithm could be.
Obviously, further applications are highly needed to see how
it may behave for solving various tough engineering optimistion problems.

\section{Conclusions}

We have successfully used the Firefly Algorithm to carry out
nonlinear design optimisation. We first validated the algorithms using
some standard test functions.  After designing some new test functions
with singularity and stochastic components, we then used the FA to solve these
unconstrained stochastic functions.
We also applied it to find
a better global solution to the pressure vessel design optimisation.
The optimisation results imply that
the Firefly Algorithm is potentially more powerful than other
existing algorithms such as particle swarm optimisation.

The convergence analysis of metaheuristic algorithms still requires some
theoretical framework. At the moment, it still lacks of a general framework for
such analysis. Fortunately, various studies started to propose
a feasible measure for comparing algorithm performance.
For example, Shilane et al (2008) suggested a framework
for evaluating statistical performance of evolutionary algorithms.
Obviously, more comparison studies are highly needed
so as to identify the strength and weakness of current
metaheuristic algorithms. Ultimately, even better optimisation
algorithms may emerge.

\end{document}